\documentclass[12pt,twoside]{amsart}
\usepackage{amsthm,amsfonts,amssymb,amscd,euscript}

\usepackage[matrix,arrow]{xy}

\author{Florin Ambro} 
\address{RIMS, Kyoto University\\
Kyoto 606-8502, Japan.}
\email{ambro@kurims.kyoto-u.ac.jp}

\newcommand{\isoto}{{\overset{\sim}{\rightarrow}}}

\newcommand{\Q}{{\mathbb Q}}
\newcommand{\Z}{{\mathbb Z}}
\newcommand{\N}{{\mathbb N}}
\newcommand{\R}{{\mathbb R}}

\newcommand{\cE}{{\mathcal E}}

\newcommand{\cO}{{\mathcal O}}
\newcommand{\cR}{{\mathcal R}}

\newcommand{\LCS}{\operatorname{LCS}}
\newcommand{\mult}{\operatorname{mult}}
\newcommand{\Proj}{\operatorname{Proj}}
\newcommand{\Supp}{\operatorname{Supp}}

\theoremstyle{plain}
\newtheorem{thm}{Theorem}[section]

\newtheorem{lem}[thm]{Lemma}
\newtheorem{cor}[thm]{Corollary}

\theoremstyle{definition}
\newtheorem{defn}[thm]{Definition}

\newtheorem{exmp}[thm]{Example}

\newtheorem{rem}[thm]{Remark}

\theoremstyle{remark}

\newenvironment{sketch}{\begin{proof}[Sketch of proof]}{\end{proof}}

\begin{document}

\bibliographystyle{amsalpha+}
\title[Log canonical centers]{Basic properties of log canonical centers}
%\date{July 16, 2006}

\begin{abstract} We present the elementary properties of log
canonical centers of log varieties.
\end{abstract}

\maketitle

%%%%%%%%%%%%%%%%%%%%%%%%%%%%%%%%%%%%%%%%%%%%%%%%%%%%%%%
%%% Document name: lccenters.tex
%%% Last modified: Wed Nov 8, 2006
%%%%%%%%%%%%%%%%%%%%%%%%%%%%%%%%%%%%%%%%%%%%%%%%%%%%%%%

%%%%%%%%%%%%%%%%%%%%%%%%%%%%%%%%%%%%%%%%
%%%%%%%%%%%%%%%%%%%%%%%%%%%%%%%%%%%%%%%%
\section{Introduction}
%%%%%%%%%%%%%%%%%%%%%%%%%%%%%%%%%%%%%%%%
%%%%%%%%%%%%%%%%%%%%%%%%%%%%%%%%%%%%%%%%

\footnotetext[1]{The author is supported by a 
21st Century COE Kyoto Mathematics Fellowship,
and a JSPS Grant-in-Aid No 17740011.
}
\footnotetext[2]{1991 Mathematics Subject Classification. 
Primary: 14B05. Secondary: 14E30
\clearpage.}

Log varieties and their log canonical centers provide a 
natural setting for inductive arguments in higher dimensional 
algebraic geometry. The prototype log variety is $(X,\sum_i b_i E_i)$, 
where $X$ is a nonsingular variety, the $E_i$'s are nonsingular 
prime divisors intersecting transversely, and $b_i\in [0,1]$. 
If $b_i<1$ for all $i$, $(X,\sum_i b_i E_i)$ has so called 
Kawamata log terminal singularities, and it has no log canonical 
centers. In general, some of the $E_i$'s will have coefficient 
one, and the connected components of their intersections
are called log canonical centers. Further, let $C$ be a log
canonical center which is a connected component 
of $\cap_{i\in J}E_i$, where $b_i=1$ for all $i\in J$.
Then a successive application of the classical adjunction 
formula gives the log adjunction formula
$$
(K_X+\sum_i b_i E_i)\vert_C=
K_C+\sum_{i\notin J}b_i E_i\vert_C.
$$
Creating log canonical centers, restricting to them by adjunction,
and lifting sections via vanishing theorems -- these are the 
three steps of a powerful technique for constructing sections 
of adjoint line bundles in characteristic zero, parallel to 
the L$^2$-methods for singular hermitian metrics in complex geometry
(see~\cite{KMM, Ein97, Kol97, A99} and ~\cite{Siu95,Dem01} for the 
algebraic and analytic side of the story, respectively). 

In this note we present the elementary properties
of log canonical centers. They are easy to see in our 
example above: its log canonical centers are nonsingular, 
finite in number, their intersections are unions of log
canonical centers, and their unions have seminormal singularities. 
Most of these properties extend to the general case of log 
varieties with log canonical 
singularities, by Shokurov~\cite{Sho92},
Koll\'ar~\cite{Kol91}, Ein-Lazarsfeld~\cite{EL93},
Kawamata~\cite{Kaw97} and~\cite{A98}, under 
some mild extra hypotheses, and~\cite{A03} in general. 
Since they seem to be obscured by the 
new terminology of quasi-log varieties in~\cite{A03}, 
we reproduce them here. 

\begin{thm}\label{basic} 
Let $(X,B)$ be a log variety with log canonical 
singularities, defined over an algebraically closed 
field of characteristic zero. Then:
\begin{itemize}
\item[(1)] $(X,B)$ has at most finitely many log canonical centers.
\item[(2)] An intersection of two log canonical centers is 
a union of log canonical centers.
\item[(3)] Any union of log canonical centers has seminormal singularities.
\item[(4)] Let $x\in \LCS(X,B)$ be a closed point. 
Then there is a unique minimal log canonical center $C_x$ 
passing through $x$, and $C_x$ is normal at $x$.
\end{itemize}
\end{thm}

In previous approaches, when there exists $0\le B^0\le B$ 
such that $(X,B^0)$ has Kawamata log terminal singularities,
Theorem~\ref{basic} follows from Kawamata-Viehweg's 
vanishing. The new idea is to use the log canonical 
version of Koll\'ar's torsion freeness, instead of vanishing.

%%%%%%%%%%%%%%%%%%%%%%%%%%%%%%%%%%%%%%%%
%%%%%%%%%%%%%%%%%%%%%%%%%%%%%%%%%%%%%%%%
\section{Log varieties}
%%%%%%%%%%%%%%%%%%%%%%%%%%%%%%%%%%%%%%%%
%%%%%%%%%%%%%%%%%%%%%%%%%%%%%%%%%%%%%%%%

A {\em log variety} $(X,B)$ is a normal variety $X$ endowed 
with an effective $\R$-Weil divisor $B$ such that $K_X+B$ is 
$\R$-Cartier. We assume that $X$ is defined over an 
algebraically closed field $k$, {\em of characteristic zero}. 
The {\em canonical divisor} $K_X$ is defined as the Weil
divisor $(\omega)$ of zeros and poles of a non-zero 
top rational differential form $\omega\in \wedge^{\dim X}
\Omega^1_X\otimes_k k(X)$ (it depends on the choice 
of $\omega$, but only up to linear equivalence). 
Now $B$ is a finite combination of 
prime divisors with real non-negative coefficients, and 
the $\R$-Cartier hypothesis means that locally 
on $X$, $K_X+B$ equals a finite sum 
$\sum_i r_i(\varphi_i)$, where $r_i\in \R$ and 
$\varphi_i\in k(X)^\times$ are non-zero rational 
functions on $X$.

Let $\mu\colon X'\to X$ be a birational morphism. If we 
use the same top rational form to define canonical 
divisors, $K_X=(\omega)$ and $K_{X'}=(\mu^*\omega)$, 
we have the {\em log pull-back formula}
$$
\mu^*(K_X+B)=K_{X'}+B_{X'}.
$$
As the notation suggests, the $\R$-Weil divisor 
$B_{X'}$ is independent of the choice of $\omega$. 
For a prime divisor $E\subset X'$, the real number
$a(E;X,B)=1-\mult_E(B_{X'})$ is called the 
{\em log discrepancy of} $(X,B)$ at $E$. This number
depends only on the valuation of $k(X)$ defined by 
$E\subset X'\stackrel{\mu}{\to} X$. We call such a
valuation {\em geometric}, and denote $c_X(E)=\mu(E)$.

\begin{defn} The log variety $(X,B)$ is said to have
\begin{itemize}
\item {\em log canonical singularities} if $a(E;X,B)\ge 0$ 
for every geometric valuation $E$ of $X$. 
\item {\em Kawamata log terminal singularities} if 
$a(E;X,B)>0$ for every geometric valuation $E$ of $X$.
\end{itemize} 
\end{defn}

The loci where $(X,B)$ has log canonical and Kawamata log 
terminal singularities, respectively, are non-empty open 
subsets of $X$. We denote their complements by 
$\LCS(X,B)$ and $(X,B)_{-\infty}$, respectively. In 
particular, $(X,B)_{-\infty}\subseteq \LCS(X,B)$,
$(X,B)_{-\infty}=\emptyset$ if and only if $(X,B)$ has
log canonical singularities, and $\LCS(X,B)=\emptyset$
if and only if $(X,B)$ has Kawamata log terminal singularities.

\begin{rem} Suppose that $(X,B)$ does not have log
canonical singularities. Then for every cycle 
$C\subseteq (X,B)_{-\infty}$ and every
positive integer $n$, there exists a geometric valuation
$E$ of $X$ such that $a(E;X,B)<-n$ and $c_X(E)=C$.
This property is behind our notation for the locus where
$(X,B)$ does not have log canonical singularities.
\end{rem}

\begin{defn}
A cycle $C\subset X$ is called a {\em log canonical center} if
$(X,B)$ has log canonical singularities at the generic point
of $C$, and there exists a geometric valuation $E$ of $X$ 
such that $a(E;X,B)=0$ and $c_X(E)=C$.
\end{defn}

\begin{exmp} Let $X$ be a nonsingular variety and
$B=\sum_{i=1}^l b_i E_i$, where $b_i\in \R_{\ge 0}$ and 
$\{E_i\}_{i=1}^l$ are nonsingular prime divisors 
intersecting transversely. Let $I_0=\{i;\ b_i=1\}$
and $I_{-\infty}=\{i;\ b_i>1\}$. Then
\begin{itemize}
\item The locus where $(X,B)$ has Kawamata log terminal
singularities is $X\setminus \cup_{i\in I_0\cup I_{-\infty}}E_i$.
Its complement $\LCS(X,B)$ has a closed subscheme structure, 
with defining ideal 
$$
\cO_X(-\sum_{i\in I_0\cup I_{-\infty}}\lfloor b_i \rfloor E_i).
$$
\item The locus where $(X,B)$ has log canonical
singularities is the open set 
$X\setminus \cup_{i\in I_{-\infty}}E_i$.
\item The log canonical centers of $(X,B)$ are the connected 
components of the intersections $\cap_{i\in J}E_i$, for
$\emptyset \ne J\subseteq I_0$.
\end{itemize}
\end{exmp}

\begin{rem} Our notion of log canonical center differs
from the standard one used in the literature. The latter 
is defined as a center $C=c_X(E)$, where $E$ is a geometric 
valuation with $a(E;X,B)\le 0$. In our case, we further
require that $(X,B)$ has log canonical singularities
at the generic point of $C$. The two notions
coincide for log varieties with log canonical singularities,
but differ otherwise. For example, consider the log variety 
$({\mathbb C}^2,H_1+2H_2)$, where $H_i:(x_i=0)$. In the standard 
literature, $H_1, H_2$, and every point of $H_2$ is a log 
canonical center. In our sense, only $H_1$ is a log canonical 
center. This seems reasonable, given that 
$K_{\mathbb C^2}+H_1+2H_2$ cannot be restricted by adjunction 
to $H_2$, or any its points.
This also shows that the finiteness in
Theorem~\ref{basic}.(1) fails in the non-log canonical case. 
\end{rem}

\begin{rem}\label{check} By Hironaka, we may choose $\mu$ so 
that $X'$ is non-singular,
and the proper transform $\mu^{-1}_*B$ and the $\mu$-exceptional
locus $\cup_{i=1}^l E_i$ is supported by a simple normal crossings
divisor. We have the following formula
$$
B_{X'}=\mu^{-1}_*B+\sum_{i=1}^l (1-a(E_i;X,B))E_i.
$$
One can see that $(X,B)_{-\infty}$ is the image in $X$ 
of the components of $B_{X'}$ with coefficients in $(1,+\infty)$,
and $\LCS(X,B)$ is the image in $X$ of the components of $B_{X'}$ 
with coefficients in $[1,+\infty)$. The lc centers of $(X,B)$ 
are the sets $\mu(S)$, where $S$ is a connected component of 
an intersection of components of $B_{X'}$ with multiplicity one, 
such that $\mu(S)\nsubseteq (X,B)_{-\infty}$. 
\end{rem}

\begin{rem} Canonical singularities were introduced by 
Reid (see~\cite{Miles87}) as the singularities 
that appear on canonical models of projective manifolds 
of general type. Likewise, log canonical singularities are the
singularities that appear on log canonical models of prototype
log varieties of general type.
It would be interesting to similarly describe {\em semi-log
canonical} singularities (see~\cite{Ale06}) and {\em quasi-log
canonical} singularities (see~\cite{A03}).

To see this in the log canonical case, let $(X,B)$ be a 
log variety with log canonical singularities. In the setting
of Remark~\ref{check}, the following formula holds:
$$
K_{X'}+\mu^{-1}_*B+\sum_{i=1}^l E_i=\mu^*(K_X+B)+
\sum_{i=1}^l a(E_i;X,B)E_i.
$$
By log canonicity, $\sum_{i=1}^l a(E_i;X,B)E_i$ is effective
and $\mu$-exceptional, so this formula becomes a relative 
Zariski decomposition for the $\mu$-big log canonical divisor 
on the left-hand side, and we obtain
$$
\bigoplus_{m\in \N}\mu_*\cO_{X'}
(m(K_{X'}+\mu^{-1}_*B+\sum_{i=1}^l E_i))
=\bigoplus_{m\in \N}\cO_X(m(K_X+B)).
$$
Further, if $B$ is rational, $(X,B)$ is recovered as
follows: the graded $\cO_X$-algebra on the left hand side
is finitely generated (since the right-hand side is), 
its $\Proj$ is $X$, and $B$ is the push forward of 
$\mu^{-1}_*B+\sum_{i=1}^l E_i$ through the natural 
map $X'\dashrightarrow \Proj$.

More generally, consider a log variety $(X',B')$ with 
log canonical singularities, and a proper morphism 
$\pi\colon X'\to S$ such that $K_{X'}+B'$ is $\pi$-big 
and rational. We expect that the graded $\cO_S$-algebra 
$\cR(X'/S,B)=\bigoplus_{m\in \N}\pi_*\cO_X(m(K_{X'}+B'))$ 
is finitely generated. If so, we obtain a natural 
birational map 
$\Phi\colon X'\dashrightarrow X:=\Proj(\cR(X'/S,B'))$,
defined over $S$, and then $(X,\Phi_* B')$ has log canonical singularities.
\end{rem}

%%%%%%%%%%%%%%%%%%%%%%%%%%%%%%%%%%%%%%%%
%%%%%%%%%%%%%%%%%%%%%%%%%%%%%%%%%%%%%%%%

\section{A torsion freeness theorem}

%%%%%%%%%%%%%%%%%%%%%%%%%%%%%%%%%%%%%%%%
%%%%%%%%%%%%%%%%%%%%%%%%%%%%%%%%%%%%%%%%
 
The main result of this section is Theorem~\ref{vt}, the 
log canonical version of Koll\'ar's torsion freeness~\cite{Kol86}.
It is a special case of~\cite[Theorem 3.2.(i)]{A03}, but we 
reproduce it here for the convenience of the reader. 
We use vanishing results of Esnault-Viehweg~\cite[Theorem 5.1]{EV92}, 
based on logarithmic de Rham complexes. Recall that the 
characteristic is zero.

\begin{thm}[\cite{EV92}, Theorem 3.2]\label{3.2}
Let $X$ be a nonsingular projective variety and $T$ a 
$\Q$-divisor such that $T\sim_\Q 0$ and 
$T-\lfloor T\rfloor=\sum_{i=1}^l \delta_i E_i$ has simple normal 
crossings support. Let $d_1,\ldots,d_l\in \Z_{\ge 0}$, and $R$ 
a reduced divisor with no common components with $\sum_{i=1}^l E_i$, 
such that $R+\sum_{i=1}^l E_i$ has simple normal 
crossings. Denote $\cE=\cO_X(-R+\lfloor T\rfloor)$.
Then the natural map of complexes
$$
\Omega^\bullet_X(\log R+\sum_{i=1}^l E_i)\otimes
\cE(-\sum_{i=1}^l d_i E_i)\to 
\Omega^\bullet_X(\log R+\sum_{i=1}^l E_i)\otimes\cE
$$ 
is a quasi-isomorphism, and the spectral sequence 
$$
E^{pq}_1=H^q(X,\Omega^p_X(\log R+\sum_{i=1}^l E_i)\otimes \cE)
\Longrightarrow 
{\mathbb H}^{p+q}(X,\Omega^\bullet_X(\log R+\sum_{i=1}^l E_i)\otimes 
\cE)
$$
degenerates at $E_1$.
\end{thm}

\begin{sketch} Let $n$ be a minimal positive integer with $nT\sim 0$.
There exists a non-zero rational function $\varphi\in k(X)^\times$ 
such that $(\varphi)=nT$. Let $\pi\colon X'\to X$ 
be the normalization of $X$ in the field $k(X)(\sqrt[n]{\varphi})$.
There exists an open subset $U\subseteq X$, with complement of
codimension at least two, such that $U'=\pi^{-1}(U)$ is nonsingular
and the restriction to $U'$ of the support $\sum_{i'}E'_{i'}$ 
of $\pi^*(\sum_i E_i)$ is also a simple normal crossings divisor.
Let $d\colon \cO_{U'}\to \Omega^1_{U'}(\log \sum_{i'}E'_{i'}\vert_{U'})$ 
be the K\"ahler differential. The homomorphism 
$\pi_*(d)$ is compatible with
action of the Galois group $\Z/n\Z=\langle \zeta\rangle$, so
its component of eigenvalue $\zeta$, denoted $\nabla$, is an 
integrable connection on $\cO_X(\lfloor T\rfloor)$ with logarithmic 
poles along $\sum_i E_i$. It is apriori defined only on $U$, but
it extends to $X$ since $\cO_X(\lfloor T\rfloor)$ is
locally free. Its residues are
$$
Res_{E_i}(\nabla)=\delta_i \cdot id\colon \cO_X(\lfloor T\rfloor)
\otimes\cO_{E_i}\to \cO_X(\lfloor T\rfloor)\otimes\cO_{E_i}.
$$
Now $\nabla$ induces an integral connection on 
$\cO_X(-R+\lfloor T\rfloor)$ with logarithmic poles along 
$R+\sum_i E_i$. Since the residue of this connection along
each $E_i$ is given by multiplication with the fractional number
$\delta_j\in (0,1)$, it follows from~\cite[Properties 2.9, 
Lemma 2.10]{EV92} that the natural map of complexes is a 
quasi-isomorphism. The last statement follows from the
degeneration of the spectral sequence associated to the 
Hodge filtration on a logarithmic de Rham 
complex (Deligne~\cite[Corollary 3.2.13]{HodgeII}),
applied to some desingularization of $X'$.
\end{sketch}

\begin{cor}[\cite{EV92}, Theorem 5.1]\label{van}
Let $L,D$ be Cartier divisors on a nonsingular projective 
variety $X$. Assume that there are nonsingular divisors 
$E_i$ intersecting transversely, and $b_i\in [0,1]$ such that:
\begin{itemize}
\item[(i)] $L\sim_\R K_X+\sum_i b_i E_i$.
\item[(ii)] $D$ is effective, supported by $\sum_{0<b_i<1}E_i$.
\end{itemize}
Then the natural maps $H^q(X,\cO_X(L))\to H^q(X,\cO_X(L+D))$ 
are injective for all $q$.
\end{cor}

\begin{proof} The assumption (i) means that there 
exist rational functions $\varphi_j\in k(X)^\times$ 
and $r_j\in \R$ such that 
$L=K_X+\sum_i b_i E_i+\sum_j r_j (\varphi_j)$. We fix
$L-K_X-\sum_{b_i\in \Q}E_i$ and regard this equality of 
divisors as a system of equations in $b_i\in \R\setminus \Q$ 
and $r_j$, in the vector space whose basis consists 
of all the prime divisors involved. The space of 
solutions is defined over $\Q$, and since a solution 
exists, a rational solution exists. Therefore we may assume that 
$b_i\in \Q$ and $L\sim_\Q K_X+\sum_i b_i E_i$. 

Denote $T=-L+K_X+\sum_i b_i E_i$. We have 
$T-\lfloor T\rfloor=\sum_{0<b_i<1}b_iE_i$ and $T\sim_\Q 0$. 
Denote $\cE=\cO_X(\lfloor -\sum_{b_i=1}E_i+T\rfloor)$
and consider the commutative diagram
\[ \xymatrix{
H^q(X,\cE (-D))  \ar[r] & H^q(X,\cE)    \\
{\mathbb H}^q(X,\Omega_X^\bullet (\log \sum_i E_i) 
\otimes \cE(-D))  \ar[u] \ar[r]^{\alpha} & 
{\mathbb H}^q(X,\Omega_X^\bullet(\log \sum_i E_i) 
\otimes \cE)\ar[u]_\beta
} \]
The first part of Theorem~\ref{3.2} gives that
$\alpha$ is an isomorphism, and the second part 
implies that $\beta$ is surjective. 
We infer that $H^q(X,\cE (-D))\to H^q(X,\cE)$ is surjective.
By Serre duality, $H^0(X,\omega_X\otimes \cE^{-1})
\to H^0(X,\omega_X\otimes \cE^{-1}(D))$ is injective.
This is the desired injective map, since 
$\omega_X\otimes \cE^{-1}=
\cO_X(\lceil L-\sum_{0<b_i<1}b_i E_i\rceil)=\cO_X(L)$.
\end{proof}

\begin{cor}\label{tako}
Let $(X,\sum_i b_i E_i)$ be a log variety such that $X$ 
is nonsingular and proper, $E_i$ are nonsingular divisors 
intersecting transversely, and $b_i\in [0,1]$ for all $i$. 
Let $L$ be a Cartier divisor on $X$, and $D$ an effective 
Cartier divisor, with the following properties:
\begin{itemize}
\item[(i)] $H=L-(K_X+\sum_i b_i E_i)$ is a semiample 
$\R$-divisor. This means that $H=\sum_i h_i H_i$, where 
$h_i\ge 0$ and $\vert H_i\vert$ are linear systems free 
of base points.
\item[(ii)] $tH-D \sim_\R D'$ for some effective $\R$-divisor 
$D'$, and $t>0$.
\item[(iii)] $\Supp(D)$ and $\Supp(D')$ do not contain 
log canonical centers of $(X,\sum_i b_i E_i)$.
\end{itemize}
Then the natural maps 
$H^q(X,\cO_X(L))\to H^q(X,\cO_X(L+D))$ are injective 
for all $q$.
\end{cor}

\begin{proof} By Hironaka, there exists a birational 
modification $\mu\colon X'\to X$ such that $X'$ is 
projective and nonsingular, and 
$B_{X'}=\mu^*(K_X+\sum_i b_i E_i)-K_{X'}$, $\mu^*D$ 
and $\mu^*D'$ are all suported by a simple normal 
crossings divisor $\sum_{i'}E'_{i'}$. Decompose
$B_{X'}=B'-A$ into the positive and negative part.

By assumption, $\mu^*D$ and $\mu^*(D')$ do not contain
components of $\lfloor B'\rfloor$. Therefore there exists
$0<\epsilon\ll 1$ such that 
$\lfloor B'+\lceil A\rceil-A+\epsilon \mu^*D+\epsilon 
\mu^*(D')\rfloor=\lfloor B'\rfloor$. Since
$\mu^*H$ is semiample, we may enlarge $\sum_{i'}E'_{i'}$ and
assume that there exists an $\R$-divisor $H'$ with the 
following properties: $H'$ is supported by 
$\sum_{i'}E'_{i'}$, it has no common components with
$B_{X'}$, $\mu^*D$ and $\mu^*(D')$, $\lfloor H'\rfloor=0$,
$H'\sim_\Q (1-\epsilon t)\mu^*H$. We obtain
$$
\lceil A \rceil+\mu^*L\sim_\Q K_{X'}+B'+\lceil A\rceil-A+
\epsilon \mu^*D+\epsilon \mu^*(D')+H'.
$$
Since the effective Cartier divisor $\mu^*D$ is supported
by the fractional part of the boundary 
$B'+\lceil A\rceil-A+\epsilon \mu^*D+\epsilon \mu^*(D')+H'$,
Corollary~\ref{van} gives the injectivity of map
$$
H^q(X',\cO_{X'}(\lceil A \rceil+\mu^*L))\to 
H^q(X',\cO_{X'}(\lceil A \rceil+\mu^*L+\mu^*D)). 
$$
On the other hand, $\mu_*\cO_{X'}(\lceil A\rceil)=\cO_X$ and 
$R^q\mu_*\cO_{X'}(\lceil A\rceil)=0$ for $q>0$. In particular,
$\mu_*\cO_{X'}(\lceil A\rceil+\mu^*L)=\cO_X(L)$ and 
$R^q\mu_*\cO_{X'}(\lceil A\rceil+\mu^*L)=0$ for $q>0$.
The Leray spectral sequence 
$$
E_2^{pq}=H^p(X,R^q\mu_*\cO_{X'}(\lceil A\rceil+\mu^*L))
\Longrightarrow H^{p+q}(X',\cO_{X'}(\lceil A\rceil+\mu^*L))
$$
degenerates, so we obtain an isomorphism
$H^q(X,\cO_X(L))\isoto H^q(X',\cO_{X'}(\lceil A\rceil+\mu^*L))$.
Similarly, we obtain an isomorphism
$H^q(X,\cO_X(L+D))\isoto H^q(X',\cO_{X'}(\lceil A\rceil+
\mu^*L+\mu^*D))$. The result follows.
\end{proof}

\begin{thm}\label{vt} Let $(X,\sum_i b_i E_i)$ be 
a log variety such that $X$ is nonsingular, $E_i$
are nonsingular divisors intersecting transversely,
and $b_i\in [0,1]$ for all $i$. Let $f\colon X\to S$
be a proper morphism, and $L$ a Cartier divisor on $X$
such that $L-K_X-\sum_i b_i E_i$ is $f$-semiample. 

Let $q\ge 0$ and $s$ a local section of $R^qf_*\cO_X(L)$,
which is zero at the generic points of
$f(X)$ and $f(C)$, for every log canonical center $C$ 
of $(X,\sum_i b_i E_i)$. Then $s=0$. 
\end{thm}

\begin{proof} We may assume that $S$ is affine and 
$f(X)=S$. Then $L-K_X-\sum_i b_i E_i$ is semiample, 
and after possibly enlarging $\sum_i E_i$, we
may assume $L\sim_\Q K_X+\sum_i b_i E_i$.

Assume by contradiction that the conclusion 
is false. Then there exists an effective very ample 
divisor $A$ on $S$ such that $f^*A$ does not contain 
any log canonical center of $(X,\sum_i b_i E_i)$, and 
the homomorphism 
$R^qf_*\cO_X(L)\to R^qf_*\cO_X(L)\otimes \cO_S(A)$
is not injective. In particular, we may compactify $X$
and $S$ and assume that the homomorphism
$$
H^0(S,R^qf_*\cO_X(L+f^*A))\to H^0(S,R^qf_*\cO_X(L+2f^*A))
$$
is not injective. 
Consider the commutative diagram of spectral sequences
\[ \xymatrix{
E^{p,q}_2= H^p(S,R^qf_*\cO_X(L+f^*A))\ar[d] &\Longrightarrow &  
H^{p+q}(X,\cO_X(L+f^*A))\ar[d]  \\
\bar{E}^{p,q}_2= H^p(S,\cO_X(L+2f^*A))  
    & \Longrightarrow  & H^{p+q}(X,\cO_X(L+2f^*A))
} \]
The map $E_2^{0,q}\to H^q(X,\cO_X(L+f^*A))$ is injective,
so we conclude that the homomorphism
$$
H^q(X,\cO_X(L+f^*A))\to H^q(X,\cO_X(L+2f^*A))
$$
is not injective. This contradicts Theorem~\ref{tako}.
\end{proof}

%%%%%%%%%%%%%%%%%%%%%%%%%%%%%%%%%%%%%%%%
%%%%%%%%%%%%%%%%%%%%%%%%%%%%%%%%%%%%%%%%

\section{Proof of Theorem~\ref{basic}}

%%%%%%%%%%%%%%%%%%%%%%%%%%%%%%%%%%%%%%%%
%%%%%%%%%%%%%%%%%%%%%%%%%%%%%%%%%%%%%%%%

\begin{rem}\label{int}
Let $X$ be a nonsingular variety, $\sum_{i\in I} E_i$ a 
simple normal crossings divisor, $I',I''$ non-empty 
subsets of $I$, and $C$ a connected component of 
$\cap_{i\in I'}E_i$. Then $C\subseteq 
\cup_{i\in I''}E_i$ if and only if $I'\cap I''\ne \emptyset$.
\end{rem}

\begin{lem}\label{ml} 
Let $(X,B)$ be a log variety with log canonical singularities, 
and $W$ a union of log canonical centers of $(X,B)$. Let 
$\mu\colon X'\to X$ be a resolution of singularities with
log pullback $\mu^*(K+B)=K_{X'}+B_{X'}$, such that 
$\mu^{-1}(W)$ is a divisor and $\mu^{-1}(W)\cup \Supp(B_{X'})$ 
has simple normal crossings.
Let $S$ be the union of prime divisors $E$ on $X'$
with $\mult_E(B_{X'})=1$ and $E\subset \mu^{-1}(W)$.
Then $\cO_W=\mu_*\cO_S$.
\end{lem}

\begin{proof} We have a unique decomposition
$
B_{X'}=S+R+\Delta-A,
$
where $S$ consists of the components of $B_{X'}$ with
coefficient equal to one and which are included in 
$\mu^{-1}(W)$, $R$ consists of the components of $B_{X'}$ 
with coefficient equal to one but which are not included 
in $\mu^{-1}(W)$, $\Delta$ is the part of $B_{X'}$ with
coefficients in $(0,1)$ and $-A$ is the negative part
of $B_{X'}$. By Remark~\ref{int} applied to $R$ and 
$\mu^{-1}(W)$, no connected component of an intersection 
of components of $R$ is mapped inside $W$. Consider 
the exact sequence
$$
\mu_*\cO_{X'}(\lceil A\rceil)\to \mu_*\cO_S(\lceil A\vert_S\rceil)
\to R^1\mu_*\cO_{X'}(\lceil A\rceil-S).
$$
By $W=\mu(S)$, $\lceil A\rceil-S=K_{X'}+R+\lceil A\rceil-A+
\Delta-\mu^*(K+B)$ and Theorem~\ref{vt}, the last map is
zero. Therefore 
$\mu_*\cO_{X'}(\lceil A\rceil)\to \mu_*\cO_S(\lceil A\rceil\vert_S)$ 
is surjective. Since $B$ is effective, we deduce that 
$A$ is $\mu$-exceptional. Therefore 
$\cO_X=\mu_*\cO_{X'}(\lceil A\rceil)$, which implies
$\mu_*\cO_S=\cO_W$.
\end{proof}

\begin{proof}[Proof of Theorem~\eqref{basic}] 
(1) Choose $\mu\colon X'\to X$ such that $X'$ is smooth 
and $B_{X'}$ has simple normal crossings support.
The log canonical centers are the images on $X$ of the components
of $B_{X'}$ with coefficient one, and their intersections.
Therefore they are finite.

(2) Let $C_1,C_2$ be two log canonical centers. By (1), it suffices 
to show that for every closed point $x\in C_1\cap C_2$, 
there exists a new log canonical center 
$x\in C_3\subset C_1\cup C_2$. 

Let $W=C_1\cup C_2$. We may choose $\mu$ so that 
the hypotheses of Lemma~\ref{ml} hold. In the notations 
of Lemma~\ref{ml}, we have
$
\mu_*\cO_S=\cO_W.
$
In particular, $S\to C_1\cup C_2$ has connected fibers.
Therefore there are prime components $E_1,E_2$ of $S$ such
that $\mu(E_i)=C_i$ and $E_1\cap E_2\cap \pi^{-1}(x)\ne 
\emptyset$. Let $Z$ be a connected component of $E_1\cap E_2$
which intersects $\pi^{-1}(x)$. Then $C_3=\mu(Z)$ is an 
log canonical center with $x\in C_3\subset C_1\cup C_2$.

(3) Let $W$ be the union of some log canonical centers.
We may choose $\mu$ so that the hypotheses of Lemma~\ref{ml} 
hold. In the notations of Lemma~\ref{ml}, we have
$
\mu_*\cO_S=\cO_W.
$
Since $S$ clearly has seminormal singularities, we infer 
by~\cite[Proposition 4.5]{A98} that $W$ has seminormal
singularities.

(4) Fix $x\in \LCS(X,B)$ and consider $(X,B)$ as a germ 
near $x$. By (1) and (3), there exists a unique log canonical 
center $x\in C$ which is minimal with respect to inclusion. 
It remains to check that $C$ is normal near $x$.
Construct $S$ as above with $\mu_*\cO_S=\cO_C$. Since 
$C$ is minimal, every connected component of an intersection 
of components of $S$ dominates $C$. 
Consider the simplicial scheme 
$(S_n=(S_0/S)^{\Delta_n}\to S)_{n\ge 0}$, 
where $S_0\to S$ is the normalization (see~\cite{HodgeIII}).
Each irreducible component of $S_n$ is nonsingular and 
mapped onto an intersection 
of components of $S$, hence it dominates $C$. Then each
$S_n\to C$ factors through the normalization of $C$.
These factorizations glue (cf.~\cite[Lemma 2.2.(ii)]{A03}),
so that $S\to C$ factors through the normalization of $C$.
But $\mu_*\cO_S=\cO_C$, so $C$ is normal.
\end{proof}

%%%%%%%%%%%%%%%%%%%%%%%%%%%%%%%%%%%%%%%
%%%%%%%%%%%%%%%%%%%%%%%%%%%%%%%%%%%%%%%

\end{document}